\pdfoutput=1
\documentclass[a4paper,11pt]{article}
\usepackage{amsmath,amscd,amsfonts,amssymb,latexsym}
\usepackage[cmtip,arrow]{xy}
\usepackage{pb-diagram,pb-xy}

\newtheorem{theorem}{Theorem}
\newtheorem{proposition}[theorem]{Proposition}

\def\C{\mathbb{C}}
\def\Q{\mathbb{Q}}

\def\F{\mathbb{F}}

\def\O{\mathbb{O}}
\def\R{\mathbb{R}}

\def\Z{\mathbb{Z}}

\def\s{\vskip10pt}

\def\dim{{\rm dim}}

\usepackage{pgf}
\def\obrazek#1{\pgfdeclareimage[height=5cm]{#1}{#1}\pgfuseimage{#1}}
\def\obrazekx#1{\pgfdeclareimage[height=3.5cm]{#1}{#1}\pgfuseimage{#1}}

\begin{document}

\author{Ma\l gorzata Mikosz\\
\small Warsaw University of Technology,\\
\small ul.~Koszykowa 75, 00-662, Warszawa, Poland\\
\small emmikosz@mini.pw.edu.pl \and
Andrzej Weber
\\
\small Department of Mathematics of Warsaw University\\
\small Banacha 2, 02-097 Warszawa, Poland\\
\small aweber@mimuw.edu.pl}

\title{\bf Triality, characteristic classes,\\  $D_4$ and $G_2$ singularities}
\date{December 2013}
\maketitle
\begin{abstract} We recall the construction of triality automorphism of $\frak{so}(8)$
given by E. Cartan and we give a matrix representation  for the real form $\frak{so}(4,4)$. We compute the induced results on the characteristic classes. Paralelly we study the triality automorphism of the singularity $D_4$ (in Arnolds classification of smooth functions) and its miniversal deformation. The similarity with Lie theory leads us to a definition of $G_2$ singularity.\end{abstract}
\s
{\small Keywords: {Triality, characteristic classes, Lie algebra $G_2$, singularities, Milnor fibration}
\s \noindent AMS Subject Classification: {
55R40, 17B25,
14J17,
32S05}}
\s
The Lie algebra $\frak{so}(8)$ is the first algebra of the series
$D_4,D_5,D_6,\dots$. It is a classical simple algebra but it is
also considered exceptional since it is the only one which admits
an automorphism of order three. The automorphism of the Lie
algebra lifts to an automorphism of the Lie group $Spin(8)$, the
universal cover of $SO(8)$. The fixed points of that automorphism
is the exceptional group $G_2$. E.~Cartan \cite{Ca} has
constructed a concrete triality automorphism for the real compact
form of $\frak{so}(8)$.  It seems that the  approach of
E.~Cartan to triality was almost forgotten. We think it is worth
to recall the original construction. In the Appendix we give a
formula for the choice of a quadratic form with signature
$(4,-4)$. It has the advantage, that the root spaces coincides
with the coordinates of the matrix.

We give formulas for the action of triality on the characteristic classes for $Spin(8)$-bundles.
It is remarkable, that the space spanned by the Euler classes of the natural representation and spin representations $S^+$, $S^-$ is a two-dimensional nontrivial representation of $\Z_3$. In other words the sum of the Euler classes is equal to zero and
$\Z_3$ permutes the Euler classes cyclically. The remaining generators of the ring of characteristic classes can be chosen to be invariant (except the case when the base field is of characteristic three) .

The Dynkin diagram $D_4$ also appears in the singularity theory.
The singularity of the type $D_4$ defined by $x^3-3xy^2$  admits
an automorphism of order three. Moreover this automorphism can be
extended to an automorphism of the parameter space of the
miniversal deformation. The action of the triality automorphism on
the functions on the parameter space is the same as the action on
the cohomology of $H^*(BSpin(8);\R)$. Taking the quotient
by the cyclic group $\Z_3$ we obtain a map germ $\C^2/\Z_3\to \C$ with a
two-dimensional  cohomology of the Milnor fiber. In the
distinguished basis of vanishing cycles the  intersections are
described by the Dynkin diagram $G_2$. We study the geometry and topology of
that singularity. Here the domain of the function is singular; it has an isolated singularity of the type $A_2$, and  taking a ramified cover, the function itself becomes the classical singularity of type $A_2$. Therefore we can say that the singularity $G_2$ is in some sense built from two singularities $A_2$. This resembles the picture of the root system of the Lie algebra $\frak{g}_2$ which contains two copies of the systems $A_2$ intertwined together.

We wish to put an emphasis on  similarities between the theory of
Lie algebras and the singularity theory. The formulas for triality
in both theories are formally the same, although they describe
objects of completely different natures.

\section{The original approach of Elie Cartan}

In 1925 E. Cartan published a paper under the title {\it Le principe de dualite et la theorie des groupes simples et semi-simples}, \cite{Ca}. In fact the main subject of the article is not  duality but a symmetry of order three. After a general introduction motivated by the duality in the matrix group $GL(n)$
$$A\mapsto (A^T)^{-1}$$
he introduces a very concrete automorphism of the Lie algebra
$\frak{so}(8)$. The construction of the transformation is the
following: Let $(a_{i,j})_{0\leq i,j \leq 7}\in \frak{so}(8)$ be
an antisymmetric matrix. For each $i\in\{1,2,\dots, 7\}$ the
quadruple
$$a_{0,i}\;,\;a_{i+1,i+5}\;,\;a_{i+4,i+6}\;,\;a_{i+2,i+3}$$ is
transformed linearly by the matrix
\begin{equation}\left(
\begin{array}{rrrr}
-\frac{1}{2} \;& -\frac{1}{2} \;& -\frac{1}{2} \;& -\frac{1}{2} \;\\ \\
 \frac{1}{2} \;& \frac{1}{2} \;& -\frac{1}{2} \;&-\frac{1}{2} \;\\ \\
 \frac{1}{2} \;& -\frac{1}{2} \;& \frac{1}{2} \;&-\frac{1}{2} \;\\ \\
 \frac{1}{2} \;& -\frac{1}{2} \;& -\frac{1}{2} \;&  \frac{1}{2}\;
\end{array}\right)\label{macierzca}\end{equation}
Here the index $i$ is understood modulo 7. The value $i=0$ plays a
special role and it is excluded from the cycle $1,2,\dots 7$. The
resulting self-map $\phi:\frak{so}(8)\to\frak{so}(8)$ is of order
three, that is $\phi\circ\phi\circ\phi=Id$ and, what is the most
important, it preserves the Lie bracket. The subalgebra
$\frak{so}(8)^\phi$ fixed by $\phi$ is the exceptional algebra of
type $\frak{g}_2$. The action of that algebra on $\R^8$
annihilates the vector $e_0$ and the trilinear form on $e_0^\perp$
$$\omega=e^*_{1 2 6}+e^*_{13 4}+e^*_{1 5 7} +e^*_{2 3 7} +
e^*_{2 4 5} +e^*_{3 5 6} +e^*_{4 67}\,,$$
where $e^*_{ijk}=e^*_i\wedge e^*_j\wedge e^*_k$. This particular form defines the imaginary part of the octonionic multiplication in $\R^8$.
Identifying $e_0$ with the unit $1\!\!1\in \O$ in octonions (Cayley numbers) and $e_0^\perp$ with $Im(\O)$  we have
$$\langle a,Im(b\cdot c)\rangle=\omega(a,b,c)$$
for $a,b,c\in Im(\O)$.
We note that $\omega$ has all coefficient equal to one, therefore this definition of octonionic multiplication  cannot agree with the most common definition via Fano plane, see e.g. \cite{Ba}. To repair the discrepancy it is enough to change the sign of the basis element $e_4$. Then the multiplication is given by the rule encoded in the picture:

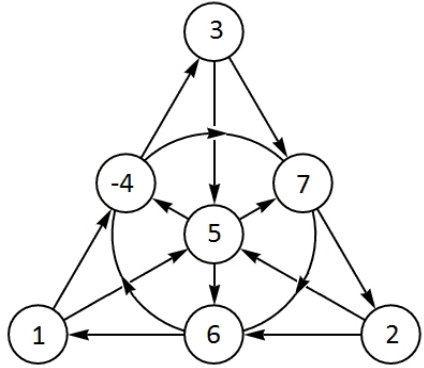
\begin{center}\obrazek{oktook}\end{center}

\noindent Moreover  $\frak{so}(8)^\phi$ is the full stabilizer of
$\omega$. It follows that the group of transformations of $\R^8$
preserving octonionic multiplication coincides with the Lie group
associated  to $\frak{so}(8)^\phi$. There are no computations in
the Cartan's paper. We refer the readers who wish to check the
formulas to \cite{MW}. A general point of view is presented and in
\cite[\S35]{KMRT}, but there the explicit form of $\phi$ is not
given. Some forms of triality is given in \cite[\S 24]{Po} or \cite[\S3.3.3]{SV}.

It is a pity that this original very simple point of view on
triality and the definition of the group $G_2$ as the fixed set
is not wide spread in the literature. Usually there are discussed
mainly the spinor representations an the nondegenerate map
$$S^+_8\times \R^8\to S^-_8$$
permuted by $\Z_3$ cyclically,
 \cite{An}, \cite[\S20.3]{Fu}, \cite{KMRT}, \cite{Vi}.
 We will construct another matrix representation of triality for $\frak{so}(4,4)$ in the Appendix.

\section{Action of triality on maximal torus}
\label{torus}
Our goal is to describe the triality in a way which does not look
like a magical trick. The approach presented  here is equivalent,
to the Cartan's work in the complex case. In of our construction
it will be clear  where the formulas come from. The triality
automorphism given below has an advantage, that the root spaces
coincide with the coordinates of the matrix and these coordinates
are permuted by $\Z_3$.

Working  with the complex coefficients we
choose a basis in $\C^8$ (as in \cite[\S19]{Fu}) in which the
quadratic form is equal to
$$x_1 x_8+x_2 x_7+x_3 x_6+x_4 x_5\,.$$ For real coefficients this
means that we deal with $SO(4,4)$. The maximal torus of $SO(4,4)$
consists of the diagonal matrices
$$diag(e^{t_1},e^{t_2},e^{t_3},e^{t_4},e^{-t_4},e^{-t_2},e^{-t_3},e^{-t_1})\,.$$
The following weights form the root system of $\frak{so}(4,4)$:
$$\pm L_i\pm L_j\quad\text{for}\;i\not=j\,,$$
where $L_i(t_1,t_2,t_3,t_4)=t_i$.
Choosing the Borel subgroup as the upper triangular matrices we obtain the Dynkin diagram of simple roots:

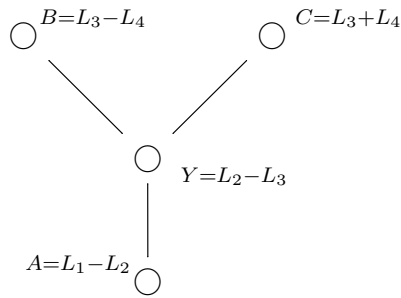
{\[
\begin{diagram}
 \node{\bigcirc}\arrow{e,l,!}{B=L_3-L_4}\arrow{se,-} \node{}\node{\bigcirc}\arrow{e,l,!}{C=L_3+L_4}\arrow{sw,-}\\
 \node{}\node{\bigcirc}\arrow{e,b,!}{~~Y=L_2-L_3}\\
 \node{}\node{\bigcirc}\arrow{w,l,!}{A=L_1-L_2}\arrow{n,-}
  \end{diagram}
\]}
The triality automorphism rotates the diagram anti-clockwise:
$$(L_1-L_2)\mapsto(L_3+L_4)\mapsto(L_3-L_4)\mapsto(L_1-L_2)$$
and fixes the root $L_2-L_3$. In the basis consisting of the weights $L_i$ the triality automorphism is given by the remarkable matrix:
\begin{equation}\left(
\begin{array}{rrrr}
\; \frac{1}{2} \;& \frac{1}{2} \;& \frac{1}{2} \;&
   -\frac{1}{2}\; \\ \\
 \frac{1}{2} \;& \frac{1}{2} \;& -\frac{1}{2} \;&
   \frac{1}{2} \;\\ \\
 \frac{1}{2} \;& -\frac{1}{2} \;& \frac{1}{2} \;&
   \frac{1}{2} \;\\ \\
 \frac{1}{2} \;& -\frac{1}{2} \;& -\frac{1}{2} \;&
   -\frac{1}{2}\;
\end{array}\right)\label{macierz}\end{equation}

\section{Action on rational cohomology of $BSO(8)$}

The rational cohomology ring of the classifying space of $SO(8)$ is well known, see e.g. \cite[Theorem 15.9]{MS}
$$H^*(BSO(n);\Q)=\Q[p_1,p_2,p_3,e]\,,$$
where $p_i$'s are the Pontryagin classes and e is the Euler class
of the tautological bundle. We identify the cohomology of the
classifying space  with the invariants of the Weyl group acting on
the polynomials on $\frak{t}$
$$H^*(BSO(n);\Q)\simeq Sym^*(\frak{t}^*)^W\,,$$
where $\frak{t}$ is the Lie algebra of the maximal torus (Cartan subalgebra),
\cite[Proposition 27.1]{Bo}. We have
$$
p_i=\sigma_i(L_1^2,L_2^2,L_3^2,L_4^2,)\quad \text{and}\quad e=L_1L_2L_3L_4\,,$$
where $\sigma_i$ is the elementary symmetric function.

\begin{theorem} Let $\phi$ be the triality automorphism acting on the cohomology of $H^*(BSO(8),\Q)=\Q[p_1,p_2,p_3,e]$. Then:
\begin{align*}
\phi(p_1)&= p_1\,,\\
\phi(p_2)&=-3 e + \frac38 p_1^2 - \frac12 p_2\,,\\
\phi(p_3)&= -\frac12 e\, p_1 +\frac1{16} p_1^3 - \frac14 p_1 p_2 + p_3
\,\\
\phi(e)&= -\frac12e - \frac1{16}p_1^2 +\frac14 p_2\,.\end{align*}

\end{theorem}

The proof is a direct computation in linear algebra. Alternatively
one can ask what are the Pontryagin and Euler classes of the spin
bundle $S^+$ associated to the universal principal bundle over
$BSpin(8)$. To know that, one has to compute the weights of the
original representation and apply the corresponding $W$-symmetric
function. For example the triality automorphism sends the Euler
class of the standard representation to the Euler class of the
spin $S^+$ representation. Its Euler class is product of weights
 {\small\begin{equation}\frac{L_1+L_2+L_3+L_4}2
 \frac{L_1+L_2-L_3-L_4}2
 \frac{L_1-L_2+L_3-L_4}2
 \frac{-L_1+L_2+L_3-L_4}2\,,\label{eulerwzor}\end{equation}}
see \cite[\S20]{Fu}.
 The result should be developed in the basis consisting of $e$,
$p_1^2$ and $p_2$.

We will chose another set of generators, which is better adapted
to homotopy theory. For any connected Lie group $G$ the cohomology
$H^*(G;\Q)$ is a Hopf algebra. By \cite{MHopf} it has to be of the
form $\Lambda P^\bullet$, where $P^\bullet$ is the graded space of
primitive elements of the Hopf algebra. (Dually, $H_*(G;\Q)$ is
isomorphic to the exterior power $\Lambda P_\bullet$, where the
primitive spaces $P_i\subset H_i(G;\Q)$ are images of homotopy
groups.) The distinguished transgression \cite[Prop.~VI, p.239]{GHV}  of the spectral sequence
$$H^p(BG;\Q)\otimes H^q(G)\Rightarrow H^{p+q}(EG)\,,$$ identifies
the primitive generators of $H^*(G;\Q)$ with the generators of
$$H^*(BG;\Q)\simeq Sym^*(P^\bullet[-1])\,,$$ where
$(P^\bullet[-1])^{i}=P^{i-1}$. The spaces
$P^{i-1}\subset H^{i}(BG;\Q)$ are preserved by automorphisms of
$G$, and therefore they form a distinguished  subspaces of
generators. For $G=SO(8)$ or $Spin(8)$ the dimensions of the
spaces $P^i$ are equal to 1 for $i=3,11$ and $\dim P^7=2$.
Therefore the action of the triality automorphism is nontrivial
only on $P^7$. We identify the space $P^7$ with a subspace of
$H^8(BSO(8);\Q)\simeq P^7\oplus S^2 P^3$.
\begin{proposition}\label{euler}
The space $P^7$ is spanned by the Euler class $e$ and its image
with respect to the triality automorphism.\end{proposition}

The orbit of the Euler class $e$ consists of $e$ and the Euler
classes of the spin representations $S^+$, $S^-$. The space
spanned by them is invariant. We will show that it is of dimension
two, in fact
$$e+e(S^+)+e(S^-)=0\,.$$
Having the expression (\ref{eulerwzor}), and the formula for the
Euler class $e(S^-)$
{\small\begin{equation}\frac{L_1+L_2+L_3-L_4}2
 \frac{L_1+L_2-L_3+L_4}2
 \frac{L_1-L_2+L_3+L_4}2
 \frac{L_1-L_2-L_3-L_4}2\,,\end{equation}}
 we check directly that
$$e(S^+)=-\frac12e - \frac1{16}p_1^2 + \frac14p_2\quad\text{and} \quad e(S^-)= -\frac12e + \frac1{16}p_1^2 -\frac14 p_2\,.$$
and the result follows.

In addition it is not hard to compute that the invariant subspace
of $H^{12}(BSO(8);\Q)$ is spanned by
\begin{equation}p_1^3\quad\text{ and}\quad p_3-\frac16 p_1 p_2\,.\label{3niez}\end{equation}

\section{Cohomology with finite coefficients}

For completeness we discuss now the cohomology with finite coefficients, although it will not be used in the remaining part of the paper.
The cohomology of $BSpin(8)$ has only 2-torsion, therefore for
$q\not=2$ the cohomology $H^*(BSpin(8);\F_q)$ is generated by the
same generators as for rational coefficients and the action of the
triality is given by the same formula. The only special issue of
for $q=3$ is the fact that $H^{12}(BSpin(8))$ as a representation
of $\Z_3$ is not semisimple. The invariant subspace spanned by the
Euler classes $e$, $e(S^+)$ and $e(S^-)$ does not admit any
invariant complement. The formula (\ref{3niez}) does not make
sense for $\F_3$.

The case $q=2$ is very different. The cohomology of $BSpin(n)$
with coefficients in $\F_2$ was computed by Quillen \cite[Theorem
6.5]{Qu}. For $n=8$ we have
 $$H^*(BSpin(8);\F_2)\simeq\F_2[w_2,w_3,w_4,w_5,w_6,w_7,w_8,w_8^+]/J\,,
$$
where $w_i$ are the Stiefel-Whitney classes of the universal
bundle and $w_8^+$ is the class of the spinor bundle. The ideal
$J$ is generated by $w_2$ and the results of Steenrod operation:
$Sq^1(w_2)$ and $[Sq^2,Sq^1](w_2)$. In general we have
\begin{align*}&Sq^1(w_2)=w_1 w_2 + w_3\,,\\ &[Sq^2,Sq^1](w_2)= w_1^3
w_2 +  w_1 w_2^2 + w_1^2 w_3 + w_2w_3 +
 w_1w_4 +  w_5\end{align*} and since here $w_1=0$ we find that
$$ H^*(BSpin(8);\F_2)\simeq\F_2[w_4,w_6,w_7,w_8,w_8^+]\,.$$

\begin{proposition} The action of the triality on $H^*(BSpin(8);\F_2)$ is the following: $\phi$ fixes $w_4$, $w_6$ and $w_7$ and
$$\phi(w_8)=w_8^+\,,\quad \phi(w_8^+)=w_8+w_8^+\,.$$
\end{proposition}

The triality has to fix $w_4$, $w_6$ and $w_7$ since the cohomology at that gradations is one dimensional. Moreover the triality permutes cyclically the representations: the natural one, $S^+$ and $S^-$, hence $\phi(w_8^+)=w_8(S^-)$. It remains to show that
$w_8(S^-)=w_8+w_8^+$. By Proposition \ref{euler} we have \begin{equation}e+e(S^+)+e(S^-)=0\label{eurow}\end{equation} in rational cohomology.
The group $H^8(BSpin(8);\F_2)$ is spanned by $w_8$, $w_8^+$ and $w_4^2$.
 The top Stiefel-Whitney classes are the reductions modulo 2 of the integral Euler classes $e$, $e(S^+)$ and $e(S^-)$. Also the class $w_4$ is the reduction modulo 2 of the integral Pontryagin class $p_1$. Therefore the relation
(\ref{eurow}) has to hold also in $\F_2$ cohomology.

We will not discuss here  the cohomology with integral
coefficients, their generators are to be found in \cite{BW}.

\section{Triality  of the singularity $D_4$}

Triality phenomenon seem  to attract recently mathematicians working on singularity theory. In a preprint \cite{IMT} the triality was related to the study of integral curves. We will discuss here only very basic and obvious appearance of triality in singularities of scalar functions.
The simple singularities of germs of holomorphic functions
$\C^n\to \C$ are indexed by the Dynkin diagrams $A_\mu$ for
$\mu\geq 1$, $D_\mu$ for $\mu\geq 4$, $E_6$, $E_7$ and $E_8$. The
Dynkin diagram $D_4$ describes the intersection form in the
homology of the Milnor fiber in the distinguished basis
corresponding to the basic vanishing cycles (defined by a choice
of paths joining the singular values of a morsification with a
regular one as in \cite[Ch. 2]{AGLV} or \cite[\S4]{Zo}).

The singularity  $D_4$  is given by the formula
$$f:\C^2\to\C\,,\qquad f(x,y)=x^3-3xy^2=x\,(x-\sqrt 3\, y)\,(x+\sqrt 3\, y)\,.$$
We have chosen the real form $D_4^-$ since the triality
does not act on the real form $D_4^+$ defined by  $x^3+3xy^2$.
We will observe how triality acts on the spaces related to that singularity. The link of the singularity consists of three circle in $S^3$ linked
with each other. They bound a surface, which is homeomorphic to the Milnor fiber:
\begin{center}\obrazek{link}\hfil\obrazek{milnor}\end{center}
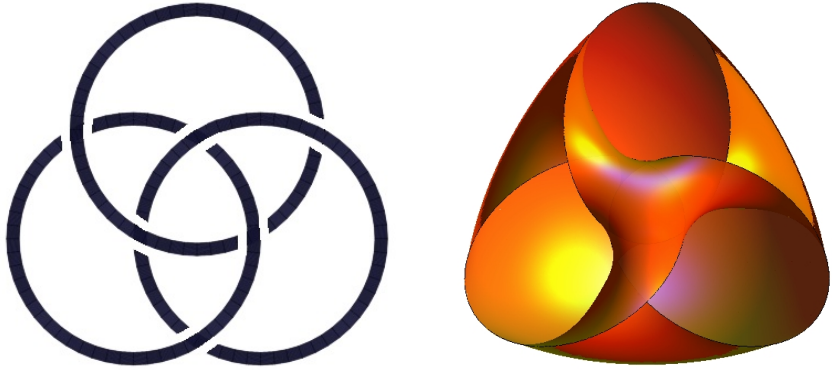
\begin{center}\textsc{ The link and the Milnor fibre\footnote{For the purpose of the picture we have taken the original  definition of  \cite{Mi} according to which the Milnor fiber is the surface in $S^3$ defined by $arg(f)=const$. Its boundary is equal to the link of the singularity.
The picture of the Milnor fiber is obtained by a stereographic projection from $S^3$ to $\R^3$.} contained in $S^3$}\end{center}

\noindent

The real part of the zero set of the function $f$ is the union of
three lines intersecting at the angle 120$^o$. The rotation of the
$(x,y)$ plane by that angle preserves the function. Denote this
rotation by $\phi_0$. The map $\phi_0$ of $\R^2$ (or $\C^2$) is
determined by the angles at which the lines intersect (up to a
cubic root of unity in the complex case). We remark, that if one
takes $x^3-xy^2 $ as the germ representing the singularity, then
the formula for $\phi_0$ and $\phi$ does not involve irrational
coefficients like $\cos(120^o)$ and $\sin(120^o)$.

The Milnor fiber $M_\varepsilon$ for $\varepsilon<<1$ is described
by the equation
$$x^3-3xy^2=\varepsilon\,,\quad |x|^2+|y|^2\leq 1\,.$$ It is an
elliptic curve with three discs removed. We can  forget the  inequality (since $f$ is homogeneous) and we identify the Milnor fibre with the plane cubic curve.  The automorphism $\phi_0$
preserves $M_\varepsilon$ permuting cyclically the removed discs.

Let
$$F:\C^2\times S\to\C$$
be a miniversal deformation of $f$, (see \cite[\S8]{AGV}, for the
definition). Here $S$, the parameter space of the deformation, it
is of dimension $\mu=4$. It can be naturally identified with
$H_1(M_\epsilon,\C)/W$, where $W$ is the Weyl group, which is
generated by the reflections in the basic vanishing cycles
corresponding to the singular points of a morsification of $f$.
The map $H_1(M_f,\C)/W\to \C^4$ is given by the period map
constructed in \cite[Theorem 1.2]{Lo},  \cite[Ch.2 \S3]{AGLV}. Let
us see how the described construction can be done equivariantly
with respect to the triality.

 First let us find a deformation which is invariant with respect
 to the rotation $\phi_0$.

\begin{proposition} The function
$$F(x,y,a,b,c,d)=x^3-3xy^2+a(x^2+y^2)+\langle(b,c),(x,y)\rangle+d$$ is a miniversal deformation of $f(x,y)=x^3-3xy^2$. It is invariant with respect to triality, provided that the action on the parameter space is given by  $$\phi(a,b,c,d)=(a,\phi_0(b,c),d)\,.$$\end{proposition}
The function $F$ is clearly invariant since $\phi_0$ preserves the
scalar product $\langle(b,c),(x,y)\rangle$. The choice of the
quadratic term is forced by the invariance condition. The
functions $1$, $x$, $y$ and $x^2+y^2$ form a basis of the local
algebra $\C[x,y]/(\frac{\partial f}{\partial x},\frac{\partial
f}{\partial y})$, as desired in the definition of miniversal
deformation.

Let us take a $\phi$-invariant morsifications of $f$
\begin{equation}f_a(x,y)=x^3-3xy^2+3a(x^2+y^2)-4a^3=(x-a)( x - \sqrt 3 y+2a) ( x - \sqrt 3 y+2a)\,.\label{morseq}\end{equation}
 We easily
compute that $f_a$ has the critical points at
$$ p_Y=(0, 0)\,,\; p_A= (-2 a, 0)\,,\;p_B=(
a , -{a}{\sqrt{3}})\,,\;p_C= (a , {a}{\sqrt{3}})\,.$$ The first
critical value is $-{4 a^3}$ (minimum), the remaining three values are equal to 0 (saddle points).
It follows that the  vanishing cycles  associated to $p_A$, $p_B$
and $p_C$ are perpendicular (see \cite[\S4, Theorem 4.26]{Zo}) they correspond to the roots $A$, $B$
and $C$. The map $\phi_0$ rotates them cyclically. Therefore
$\phi_0$ realizes the triality automorphism of the root system
 $D_4$. Topologically the system of vanishing cycles is homeomorphic to the following configuration:
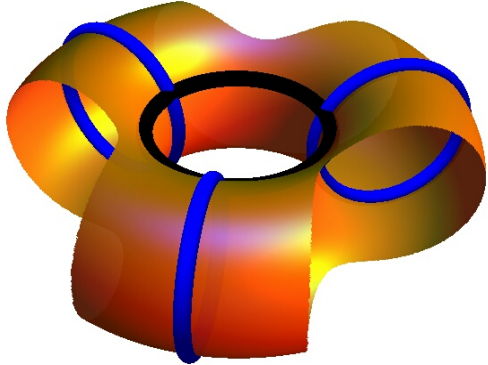
\begin{center}\obrazek{D4abs}\end{center}
\noindent
Note that the action of $\phi$ coincides with the action on the distinguished generators of $H^*(BSpin(8);\R)$ found in
Proposition \ref{euler}. The singularity $D_4$ is homogeneous, and the coefficients of the miniversal deformation can be given the gradations
$$\deg(a)=1,\quad \deg(b)=\deg(c)=2,\quad\deg(d)=3\,,$$
so that the whole function $F$ is quasihomogeneous. We obtain

\begin{theorem} Let  $a,b,c,d$ be the coefficients of the miniversal deformation of the singularity $D_4^-$. There exists an isomorphism
$$H^*(BSO(8),\R)\simeq \R[a,b,c,d]\,, $$ which preserves the action of the triality automorphism and agrees with the gradation after dividing by 4 the degrees in $H^*(BSO(8),\R)$.
\end{theorem}
If fact this isomorphism is the composition of the isomorphisms $$H^*(BSO(8),\R)\simeq Sym^*(\frak{t}^*)^W\simeq Sym^*(H^1(M_\epsilon))^W\simeq \C[H_1(M_\epsilon)/W]$$
and the period mapping.
 The last step is the most involving since according to the assumption of \cite[Theorem 1.2]{Lo} demands enlarging the dimension of the domain of $f$. We just settle for an abstract isomorphism.

The bifurcation set of the singularity $D_4^-$  is defined by
$$\{(a,b,c)\in\R^3\;|\;\exists \,(x, y)\in\R^2 \;grad(F)(x,y)=0,
\;Hessian(F)(x,y)=0\}$$  in the reduced
parameter space $d=0$. It has clearly a $\Z_3$ symmetry.
\begin{center}\obrazek{piramida}\end{center}
\begin{center}\textsc{The pyramid}\end{center}

\section{Singularity $G_2$?}
The Dynkin diagram $G_2$ does not appear in the original Arnold's
classification of simple singularities as well as the series
$B_\mu$, $C_\mu$ and the exceptional $F_4$. The other diagrams
appear in \cite{Ar} (see also references in \cite{GH}), while
$G_2$ is only mentioned in remark at the end of $\S9$. The series
$B_\mu$ and $C_\mu$ arise as diagrams for the singularities with
boundary condition. We will show how the diagram $G_2$ appears for
singularities with $\Z_3$ symmetry.

It is worth to continue
the analogy with the world of Lie algebras. Here the situation is dual, instead of taking the fixed points we divide by the $\Z/3$ action. The function $f$
factors to the quotient $\bar f: \C^2/\Z_3\to \C$. The quotient
space is not smooth, but it has mild singularities, an isolated du
Val singularity of the type $A_2$. The new Milnor fiber $\overline
M_\varepsilon$ is the quotient of the original Milnor fiber
$M_\varepsilon$. The quotient map is an unbranched cover.
Therefore the $\overline M_\varepsilon$ is homeomorphic to an
elliptic curve with one disc removed. The cohomology
$H^1(\overline M_\varepsilon;\Q)$ is generated  by two vanishing
cycles corresponding to the singular values of the invariant morsification
(\ref{morseq}). The vanishing cycle corresponding to the value $0$
is the usual one. The value $-4 a^3$ corresponds to the vanishing
cycle shrieked to the singular point of the domain. It is
reasonable to treat the quotient space as $\C^2/\Z_3$ as a stack.
Here it simply means  that
we consider the singularity $D_4$ together with the
$\Z_3$-symmetry, as it was done in \cite{Go} for unitary
reflection groups. We compute the
intersection number of a pair of cycles taking their inverse images
in the cover and dividing the result by the order of the cover.  \begin{center}\obrazekx{G2D4abs}\end{center}
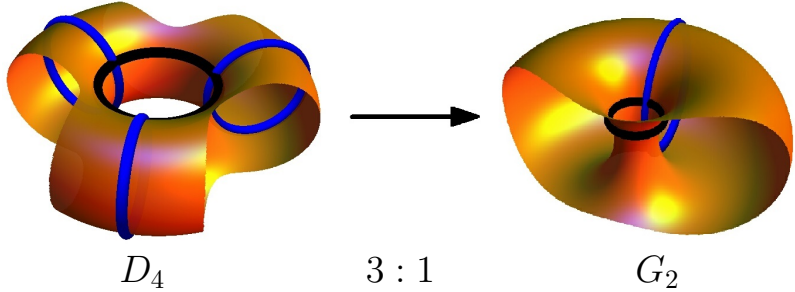
\begin{center}\vskip-20pt
\Large$D_4$\hfil 3\;:\;1\hfil$G_2$\end{center}
\noindent
To see what
is the associated diagram of that singularity we pass to an odd
dimension, adding a nondegenerate quadratic form, e.g.
$$\tilde f(x,y,z)=x^3-3xy^2+z^2\,.$$
Then the self-intersection of the cycle corresponding to $[p_A]$
is $-2$, and the self intersection of the cycle corresponding to
$[p_Y]$ is equal to $-\frac23$. The intersection diagram is
exactly $G_2$:
$$_{[p_A]=[p_B]=[p_C]}\;\bigcirc
\equiv\!\!\!\equiv\!\!\!\equiv\!\!\!\!\big>\!\!\!\equiv\!\!\!\equiv
\bigcirc\;_{[p_Y]}.$$
We will give more precise description of the Milnor fiber of $\bar f$. It is contained in the quotient space $\C^2/\Z_3$ which has the
ring of algebraic functions isomorphic to
$$\C[u,v,w]/(u^2+v^2-w^3)\,,$$
by setting $u=x^3-3x y^2$, $v=y^3 -3 x^2 y$ and $w=x^2+y^2$.
The function $\bar f$ is equal to $u$. We see that $\overline {M}_\epsilon$ is described in $\C^3$ by the equations
$$\left\{\begin{matrix}u^2+v^2&=&w^3\cr u&=&\varepsilon\,.\end{matrix}\right.$$
Hence it is isomorphic to the plane cubic $$v^2=w^3-\varepsilon^2\,.$$
All the fibers are isomorphic as algebraic curves.
Topologically, the classical Milnor fibration
$$arg(f):S^3\setminus f^{-1}(0)\to S^1$$
descends to a fibration from the complement of a circle in the lens space
$$arg(f):L(3;1)\setminus S^1\to S^1$$
with the fiber homeomorphic to the topological 2-torus with one point removed.

Let $g(s,t)=s^3-t^2$ be a germ of $A_2$ singularity.
Note that we have a commutative diagram
$$\begin{matrix}\C^2/\Z_3&\stackrel{(w, v)}\longrightarrow &\C^2\\
_{\bar f}\downarrow\phantom{^{\bar f}}&&
\phantom{^ g\,}\downarrow\,_g \\
\C&\stackrel{z\mapsto z^2}\longrightarrow &\C\,.\end{matrix}
$$
Therefore the germ of $\bar f$ is induced from $g$ by a double covering of $\C$. It is not coincidence. The root system of the Lie algebra $\frak{g}_2$ contains the root system of the algebra of the type $A_2$.

The miniversal deformation of $G_2$ germ is understood as a slice in the space of jets of invariant function, which is transverse to the orbit of the function $f$. It is clearly equal to
$$F[x,y,a,0,0,d]=u+a\, w+d\,.$$
The associated family of fibers
$$\left\{\begin{matrix}u^2+v^2&=&w^3\cr u+a\,w+d&=&0\,.\end{matrix}\right.$$
is a nontrivial family of plane cubic curves
$$v^2=w^3-(a\,w+d)^2\,.$$

We conclude with the remark that an analogous  construction can be applied to the singularity of $n$ lines intersecting at one point
$$f(x,y)=\prod_{k=1}^n(\cos({k\pi}/n) x+\sin({k\pi}/n)y ),$$ which is symmetric with respect to the group $\Z_n$ acting by the rotations with
the angles $\frac{2k\pi}n$. The resulting quotient Milnor fiber is isomorphic to a complex curve of genus $\frac{n-1}2$ with one disc removed, when $n$ is odd. For $n$ even the quotient Milnor fiber is a complex curve of genus $\frac{n-2}2$ with two discs removed. That construction for $n>3$ seem not to have a counterpart in the realm of Lie algebras.
\section{Appendix}
\subsection{Triality in $\frak{so}(4,4)$}
In Section \S\ref{torus} we have given a formula for triality automorphism acting on  the dual $\frak{t}^*$ of the Cartan subalgebra of $\frak{so}(4,4)$.
It does not preserve the lattice
corresponding to the group $SO(4,4)$ but it preserves the lattice
spanned by $L_i$'s and $\frac12(L_1+L_2+L_3+L_4)$, which
corresponds to  $Spin(4,4)$, the cover of $SO(4,4)$. We list below
the set of positive roots:
$$
\begin{array}{cccccccc}
  & L_{1}-L_{2} & L_{1}-L_{3} & L_{1}-L_{4} & L_{1}+L_{4} & \boxed{L_{1}+L_{3}} &\boxed{ L_{1}+L_{2} }&  \\
 &  & \boxed{L_{2}-L_{3}} & L_{2}-L_{4} & L_{2}+L_{4} & L_{2}+L_{3} &
   &  \\
  &  &  & L_{3}-L_{4} &L_{3}+L_{4} &  &  &
\end{array}
$$
The boxed roots are fixed by triality. It can be easily seen
when we express the roots in the basis of simple roots:
{\def\pp{\phantom{ABC2Y}}$$
 \begin{array}{cccccc}
  A & AY & ABY & ACY &\boxed{ABCY} & \boxed{ABC2Y}\\
   & \boxed{Y} & BY & CY & BCY &   \\
   &  & B & C &  &    \\
 \pp &\pp  &\pp  & \pp & \pp &\pp
\end{array}
$$}
\vskip-20pt
\noindent Here for example $ ABC2Y$ denotes the root $A+B+C+2Y$.
Dividing roots into orbits of the triality automorphism we see that we have
the fixed roots
$Y$, $ABCY$ and $ABC2Y$.
Three free orbits are generated by
$A$, $AY$ and $ABY$.

From  general theory it follows that the triality automorphism
of weights lifts to a self-map of the Lie algebra
$\frak{so}(4,4)$. But a priori it is not clear that one can find
such a lift of order three. Not every lift satisfies
$\phi\circ\phi\circ\phi=Id$. The choice of signs is not obvious
and demands a careful check. All the calculations  can be found in
\cite{MW}. The elements of $\frak{so}(4,4)$ for our quadratic form
defined by the matrix with 1's on the antidiagonal are the
matrices $\left(m_{ij}\right)_{1\leq i,j\leq 8}$ which are
antisymmetric with respect to the reflection in the antidiagonal.
Such a matrix is transformed by the triality automorphism to the
following one
$$\left(
\begin{array}{cccccccc}
\bullet & m_{34} & -m_{24} &
   m_{26} & m_{14} & \boxed{m_{16}} & \boxed{m_{17}} & 0 \\
 m_{43} & \bullet & \boxed{m_{23}} &
   m_{25} & -m_{13} & m_{15} & 0 & \boxed{-m_{17}} \\
 -m_{42} & \boxed{m_{32}} & \bullet &
   m_{35} & m_{12} & 0 & -m_{15} & \boxed{-m_{16}} \\
 m_{62} & m_{52} & m_{53} &\bullet & 0 & -m_{12} & m_{13} & -m_{14}
   \\
 m_{41} & -m_{31} & m_{21} & 0 &\bullet & -m_{35} & -m_{25} & -m_{26} \\
 \boxed{m_{61}} & m_{51} & 0 & -m_{21} & -m_{53} & \bullet & \boxed{-m_{23}} & m_{24} \\
 \boxed{m_{71}} & 0 & -m_{51} & m_{31} & -m_{52} & \boxed{-m_{32}} & \bullet & -m_{34} \\
 0 & \boxed{-m_{71}} & \boxed{-m_{61}} & -m_{41} & -m_{62} & m_{42} & -m_{43} &
   \bullet
\end{array}
\right)$$ The upper half of the diagonal
$(t_1,t_2,t_3,t_4)=(m_{11},m_{22},m_{33},m_{44})$ is transformed
by the matrix (\ref{macierz}) to
\begin{align*}
 \frac12 (m_{11} + m_{22} + m_{33} -& m_{44})\\
 \frac12 (m_{11} + m_{22} - &m_{33} + m_{44})\\
 \frac12 (m_{11} -& m_{22} + m_{33} + m_{44})\\
 \frac12 (&m_{11} - m_{22} - m_{33} - m_{44})\end{align*}
We remark that this is an unique automorphism with real
coefficients extending the self-map of the maximal torus. We would
like to stress, that both: the Cartan construction of triality for
$\frak{so}(8)$ and the triality for $\frak{so}(4,4)$ presented
here works for any ring in which 2 is invertible. On the other
hand we  note that the remaining real forms of the orthogonal
algebra $\frak{so}(k,k-8)$ for $k=1,2,3,5,6,7$ do not admit any
triality automorphism with real coefficients. Already on the level
of $\frak t^*$ we obtain the rotation matrices with imaginary
coefficients.

\subsection{Noncompact version of $G_2$}
The algebra fixed by $\phi$ consist of matrices of the form
$$\left(
\begin{array}{cccccccc}
 {t_2}+{t_3} & {a_1} & -{a_2} & {a_3} &
   {a_3} & {a_4} & {a_5} & 0 \\
 {a_6} & {t_2} & {a_7} & {a_2} & {a_2} &
   {a_3} & 0 & -{a_5} \\
 -{a_8} & {a_9} & {t_3} & {a_1} & {a_1}
   & 0 & -{a_3} & -{a_4} \\
 {a_{10}} & {a_8} & {a_6} & 0 & 0 & -{a_1} &
   -{a_2} & -{a_3} \\
 {a_{10}} & {a_8} & {a_6} & 0 & 0 & -{a_1} &
   -{a_2} & -{a_3} \\
 {a_{11}} & {a_{10}} & 0 & -{a_6} & -{a_6} &
   -{t_3} & -{a_7} & {a_2} \\
 {a_{12}} & 0 & -{a_{10}} & -{a_8} & -{a_8} &
   -{a_9} & -{t_2} & -{a_1} \\
 0 & -{a_{12}} & -{a_{11}} & -{a_{10}} & -{a_{10}} &
   {a_8} & -{a_6} & -{t_2}-{t_3}
\end{array}
\right)$$ Below we formulate a fundamental fact (known to the specialists) to which we could
not find a reference except  the  paper of \cite{Ca} for $\frak{so}(8)$ form or
without a explicit matrix realization in \cite[\S 24]{Po}, \cite[\S3.3.3]{SV}, \cite[Prop
35.9]{KMRT}.

\begin{theorem}The fixed points of the triality automorphism is a Lie algebra of the type $\frak{g}_2$.\end{theorem}

Indeed, the fixed elements of $\frak{t}^*$ is spanned by the simple roots $Y=L_2-L_3$ (the longer root) and $\frac13(A+B+C)=\frac13(L_1-L_2+2L_3)$ (the shorter root)
$$_{\frac13(A+B+C)}\;\bigcirc
\equiv\!\!\!\equiv\!\!\!\equiv\!\!\!\!\big<\!\!\!\equiv\!\!\!\equiv\bigcirc\;_Y\,.$$
The shorter root as the functional on $\frak{t}^\phi$ is equal to $A$ or $B$ or $C$, but we represent it as an invariant element $\frac13(A+B+C)\in(\frak{t}^*)^\phi$.
The  positive root spaces of $\frak{so}(4,4)^\phi$ are the
following:\begin{itemize} \item the eigenspaces associated to the
longer roots $$\frak{so}(4,4)_{ Y},\quad\frak{so}(4,4)_{
ABCY},\quad\frak{so}(4,4)_{ ABC2Y}$$
\item the diagonal subspaces associated to the shorter roots $\frac13(A+B+C)$, $\frac13(AY+BY+CY)$, $\frac13(ABY+ACY+BCY)$
$$\big(\frak{so}(4,4)_{ A}\oplus\frak{so}(4,4)_{ B}\oplus\frak{so}(4,4)_{ C}\big)^\phi
$$ $$
\big(\frak{so}(4,4)_{ AY}\oplus\frak{so}(4,4)_{ BY}\oplus\frak{so}(4,4)_{ CY}\big)^\phi,
$$ $$
\big(\frak{so}(4,4)_{ ABY}\oplus\frak{so}(4,4)_{ ACY}\oplus\frak{so}(4,4)_{ BCY}\big)^\phi.$$
\end{itemize}

\subsection{The twin-brother of octonion algebra}
We repeat the construction of the $\frak g_2$ algebra and the
octonions from Section~1. The Lie group corresponding to
$\frak{so}(4,4)^\phi$ is the noncompact form of $G_2$ since its
maximal torus is $(\R^*)^2$. The Lie algebra $\frak{so}(4,4)^\phi$
annihilates the vector $v_0=e_4-e_5$ and its action restricted to
$v_0^\perp$ preserves the form
$$\widetilde{\omega}=2e^*_{1 6 7} + 2 e^*_{2 3 8} +
 e^*_{2 (4+5) 7}-e^*_{1 (4+5) 8} + e^*_{3 (4+5) 6}
\,,$$ where $e^*_{i(4+5)j}$ means $e^*_i\wedge (e^*_4+e^*_5)\wedge
e^*_j$. The associated algebra is the algebra of pseudo-Cayley
numbers (also called split octonions). The form $\widetilde{\omega}$ is of the type 5 according
to the classification of real multisymplectic forms in dimension 7
in \cite{BV}.

\end{document}